\documentclass[12pt]{article}
\usepackage{amsfonts}
\usepackage{amssymb}
\usepackage{mathrsfs}
\usepackage{amsmath}
\usepackage{amsthm}
\usepackage{cite}

  \date{ }
     \title{{\bf Generalized  Derivations of Hom-Lie Superalgebras}
   \thanks{Supported by  NNSF of China (No. 11171055),  NSF of  Jilin province (No.201115006), Scientific
    Research Foundation for Returned Scholars Ministry of Education of China.  Corresponding author (L. Chen):
     chenly640@nenu.edu.cn}}
   \author{Jia Zhou, Liangyun Chen, Yao Ma
     \vspace{0.3cm} \\School of Information Technology, Jilin Agricultural University,\\ Changchun, 130024, CHINA \\School of Mathematics and Statistics, Northeast Normal University,\\Changchun, 130024, CHINA}
 
\begin{document}
 \maketitle
 \begin{center}{\bf Abstract}\end{center}


  In this paper, we give some basic properties
of the generalized derivation algebra ${\rm GDer}(L)$ of a Hom-Lie
superalgebra $L$. In particular, we prove that ${\rm GDer}(L) = {\rm
QDer}(L) + {\rm QC}(L)$, the sum of the quasiderivation algebra and
the quasicentroid. We also prove that ${\rm QDer}(L)$ can be
embedded as derivations in a larger Hom-Lie superalgebra.
   \vspace{0.3cm}

    \vspace{0.3cm}


   \vspace{0.3cm}
  \noindent {\bf \S 0\quad Introduction}
 \vspace{0.5cm}

  Hom-Lie algebras are a generalization of Lie algebras, where the
  classical Jacobi identity is twisted by a linear map. In the particular case that the twisted map is the identity map, Hom-Lie algebras become Lie algebras. The notion of
  Hom-Lie algebras was introduced by Hartwig, Larsson and
  Silvestrov to describe the structures on certain deformations of
  the Witt algebras and the Virasoro algebras \cite{4}. Hom-Lie algebras
  are also related to deformed vector fields, the various versions
  of the Yang-Baxter equations, braid group representations, and
  quantum groups \cite{4,5,6}. Recently, Hom-Lie superalgebras were studied in
  \cite{8,9,Liu&Chen&Ma}. More applications of the Hom-Lie algebras,
  Hom-algebras and Hom-Lie superalgebras can be found in \cite{10,11,12}.

   As is well known, derivation and generalized derivation algebra are very
  important subjects both in the research of Lie algebras and Lie
  superalgebras. In the study of Levi factors in derivation algebra
  of nilpotent Lie algebras, the generalized derivations,
  quasiderivations, centroids, and quasicentroids play key roles \cite{13}.
    Melville dealt particularly with the
  centroids of nilpotent Lie algebras \cite{14}. The most important and
systematic research on the generalized derivation algebra of a Lie
algebra was due to Leger and Luks \cite{15}. In \cite{13}, some nice properties of the generalized derivation algebra and their subalgebras, for example, of the quasiderivation algebra and of the centroid have been obtained. In particular,
they investigated the structure of the generalized derivation algebra and characterized the Lie algebras satisfying certain
conditions. Meanwhile, they also pointed that there exist some
connections between quasiderivations and cohomology of Lie algebras.

 For the generalized derivation algebra of general
non-associative algebras, the readers will be referred to \cite{14,21,22,27}.

The purpose of this paper is to generalize some beautiful results in
\cite{13,27} to the generalized derivation algebra of a Hom-Lie
superalgebra. In this paper, we mainly study the derivation algebra
${\rm Der}(L)$, the center derivation algebra ZDer($L$), the
quasiderivation  algebra QDer($L$),  and the generalized derivation
algebra GDer($L$) of a Hom-Lie superalgebra $L$.

We proceed as follows. Firstly we recall some basic definitions and
propositions which will be used in what follows. Then we give some
basic properties of the generalized derivation algebra and their
Hom-subalgebras, show that the quasicentroid of a Hom-Lie
superalgebra is also a Hom-Lie superalgebra if only and if it is a
Hom-associative superalgebra. Finally we prove that the
quasiderivations of $L$ can be embedded as derivations in a larger
Hom-Lie superalgebra  and obtain a direct sum decomposition of
${\rm Der}(L)$ when the annihilator of $L$ is equal to zero.

\vspace{0.3cm}
  \noindent {\bf \S 1\quad Preliminaries}
 \vspace{0.5cm}

 Throughout this paper $\bf K$ is a field of characteristic zero
 and $\rm Z_{2}=\{\bar{0},\bar{1}\}$. A vector space $V$ is said
 to be $\rm Z_{2}$-graded if $V=V_{\bar{0}}\oplus V_{\bar{1}}$. An element $x\in
 V_{\gamma}~(\gamma\in\rm Z_{2})$ is said to be homogeneous of degree $\gamma$. For
 simplicity the degree of an element $x$ is denoted by $x$. Let $V$
 and $W$ be two $\rm Z_{2}$-graded vectors spaces. A linear map
 $f:V\rightarrow W$ is said to be homogeneous of degree
 $\xi\in\rm Z_{2},$ if $f(x)$ is homogeneous of degree $\gamma +\xi$
for all the element $x\in V_{\gamma}$. The set of all such maps
 is denoted by ${Hom(V,W)}_{\xi}$. It is a subspace of $Hom(V,W)$,
 the vector space of all linear maps from $V$ into $W$. If in
 addition, $f$ is homogeneous of degree $\bar{0}$, i.e.$f(V_{\gamma})\subseteq
 V^{'}_{\gamma}$ holds for any $\gamma\in\rm Z_{2}$, then $f$ is said to
 be even.
\vspace{0.3cm}

 \noindent{\bf Definition 1.1}\, \cite{1} \quad{\it A Hom-Lie superalgebra is a triple $(L,[\cdot,\cdot],\alpha)$ consisting of a $\rm Z_2$-graded vector space $L=L_{\bar{0}}\oplus
L_{\bar{1}},$ an even bilinear map $[\cdot,\cdot]:L\times L\rightarrow L$
(i.e. $[L_{\theta},L_{\mu}]\subseteq L_{\theta+\mu}$) and an even
linear map $\alpha:L\rightarrow
  L$ such that for homogeneous elements $x,y,z\in L$ we have

 $(1)$\, $[x,y]=-(-1)^{xy}[y,x], $

 $(2)$\, $(-1)^{zx}[\alpha (x),[y,z]]+(-1)^{xy}[\alpha (y),[z,x]]+(-1)^{yz}[\alpha (z),[x,y]]=0.$

\noindent In particular, if $\alpha$ preserves the bracket,
(i.e. $\alpha[x,y]=[\alpha x,\alpha y],~\forall x,y\in L$), then we
 call $(L,[\cdot,\cdot],\alpha)$ a multiplicative Hom-Lie
 superalgebra.}
\vspace{0.3cm}

 \noindent{\bf Definition 1.2}\,  \cite{1}\quad  {\it Let
$(L,[\cdot,\cdot],\alpha)$ be a Hom-Lie superalgebra and define the
following subvector space $\mho$ of ${\rm End}(L)$ consisting of even
linear maps $u$ on $L$ as follows:
 $$\mho=\{u\in{\rm End}(L) |~u\alpha=\alpha u \}$$and $$\sigma:\mho\rightarrow \mho;~\sigma(u)=\alpha
 u.$$ Then $\mho$ is a Hom-Lie superalgebra over $\bf K$ with the bracket
$$[D_{\theta},D_{\mu}]=D_{\theta}D_{\mu}-(-1)^{\theta\mu}D_{\mu}D_{\theta}$$ for all $D_{\theta},D_{\mu}\in {\rm hg}(\mho).$}\vspace{0.3cm}

\noindent{\bf Definition 1.3}\,  \cite{1}\quad Let {\it
$(L,[\cdot,\cdot],\alpha)$ be a multiplicative Hom-Lie superalgebra. A
homogeneous bilinear map $D:L\rightarrow L$ of degree $d$ is said to
be an $\alpha^{k}$-derivation of $L$, where $k\in\bf N$, if it
satisfies
$$D\alpha=\alpha D,$$
$$[D(x),\alpha^{k}(y)]+(-1)^{dx}[\alpha^{k}(x),D(y)]=D([x,y]),$$
$\forall x\in {\rm hg}(L),~y\in L.$}

We denote the set of all $\alpha^{k}$-derivations by ${\rm
Der}_{\alpha^{k}}(L)$, then ${\rm Der}(L):=\bigoplus_{k\ge0}{\rm
Der}_{\alpha^{k}}(L)$ provided with the super-commutator and the
following even map
$$\tilde \alpha:\rm Der(L)\rightarrow Der(L);~~\tilde \alpha(D)=D\alpha$$
is a Hom-subalgebra of $\mho$ and is called the derivation algebra
of $L$.\vspace{0.3cm}

 \noindent{\bf Definition 1.4}\, \cite{2}\quad {\it An endomorphism $D\in {\rm Der}_{\theta}(L)$ is said to be a homogeneous generalized
 $\alpha^{k}$-derivation of degree $\theta$ of $L$, if there exist two
 endomorphisms
$D',D''\in {\rm End}_{\theta}(L)$ such that
$$D\alpha=\alpha D=0,~D\alpha^{'}=\alpha^{'} D=0,~D\alpha^{''}=\alpha^{''} D=0$$
$$[D(x),\alpha^{k}(y)]+(-1)^{\theta x}[\alpha^{k}(x),D'(y)]=D''([x,y]),\eqno(1.1)$$
for all $x\in {\rm hg}(L), y\in L.$} \vspace{0.3cm}

 \noindent{\bf Definition 1.5}\, \cite{2}\quad {\it An endomorphism $D\in {\rm
End}_{\theta}(L)$ is said to be a homogeneous $\alpha^{k}$-quasiderivation of
degree $\theta$, if there exists an endomorphism $D'\in {\rm
End}_{\theta}(L)$ such that
$$D\alpha=\alpha D=0,~D\alpha^{'}=\alpha^{'} D=0,$$
$$[D(x),\alpha^{k}(y)]+(-1)^{\theta x}[\alpha^{k}(x),D(y)]=D'([x,y]),\eqno(1.2)$$ for all $x\in {\rm hg}(L), y\in
L.$}

Let ${\rm GDer}_{\alpha^{k}}(L)$ and ${\rm QDer}_{\alpha^{k}}(L)$ be
the sets of homogeneous generalized $\alpha^{k}$-derivations and of
homogeneous $\alpha^{k}$-quasiderivations, respectively. That is
$${\rm GDer}(L):=\bigoplus_{k\ge 0}{\rm GDer}_{\alpha^{k}}(L), \quad{\rm
QDer}(L):=\bigoplus_{k\ge 0}{\rm QDer}_{\alpha^{k}}(L).$$
 It is easy to verify that both ${\rm GDer}(L)$ and
${\rm QDer}(L)$ are Hom-subalgebras of  $\mho$ (see Proposition
$2.1$).\vspace{0.3cm}

 \noindent{\bf Definition 1.6}\, \cite{2}\quad {\it If ${\rm C}(L):=\bigoplus_{k\ge 0}{\rm
C}_{\alpha^{k}}(L),$ with ${\rm C}_{\alpha^{k}}(L)$ consisting of
$D\in hg({\rm End}(L))$ of degree $d$ satisfying$$D\alpha=\alpha
D,$$
$$[D(x),\alpha^{k}(y)]=(-1)^{dx}[\alpha^{k}(x),D(y)]=D([x,y]),$$
for all $x\in {\rm hg}(L),y\in L,$ then ${\rm C}(L)$ is called an
$\alpha^{k}$-centroid of $L$.}\vspace{0.3cm}

 \noindent{\bf Definition 1.7}\, \cite{2}\quad {\it If
${\rm QC}(L):=\bigoplus_{k\ge 0}{\rm QC}_{\alpha^{k}}(L)$ with ${\rm
QC}_{\alpha^{k}}(L)$ consisting of $D\in hg({\rm End}(L))$ of degree
$d$ such that
$$[D(x),\alpha^{k}(y)]=(-1)^{dx}[\alpha^{k}(x),D(y)],$$ for all
$x\in {\rm hg}(L),y\in L, $ then ${\rm QC}(L)$ is called an
$\alpha^{k}$-quasicentroid of $L$.}\vspace{0.3cm}

Define ${\rm ZDer}(L):=\bigoplus_{k\ge 0}{\rm Der}_{\alpha^{k}}(L)$,
where ${\rm Der}_{\alpha^{k}}(L)$  consists of $D\in hg({\rm
End}(L))$ such that $$[D(x),\alpha^{k}(y)]=D([x,y])=0,$$ for all
$x\in {\rm hg}(L),~y\in L.$

It is easy to verify that $${\rm ZDer}(L)\subseteq {\rm
Der}(L)\subseteq {\rm QDer}(L)\subseteq {\rm GDer}(L)\subseteq {\rm
Pl}(L).$$
$${\rm C}(L)\subseteq {\rm QC}(L)\subseteq {\rm QDer}(L).$$

\noindent{\bf Definition 1.8}\, \cite{1}\quad {\it Let
$(L,[\cdot,\cdot],\alpha)$ be a multiplicative Hom-Lie superalgebra. If
${\rm Z}(L):=\bigoplus_{\theta\in \Gamma}{\rm Z}_{\theta}(L),$ with
${\rm Z}_{\theta}(L)=\{x\in L_{\theta}| [x,y]=0, \forall x\in {\rm
hg}(L),y\in L\},$ then ${\rm Z}(L)$ is called the center of
$L$.}\vspace{0.3cm}

     \noindent{\bf \S 2\quad Generalized derivation algebra and their  Hom-subalgebras}
 \vspace{0.5cm}

  First, we  give some basic properties
of center derivation algebra, quasiderivation algebra and the
generalized derivation algebra of a Hom-Lie superalgebra.

     \noindent{\bf  Proposition 2.1}\quad {\it Let
$(L,[\cdot,\cdot],\alpha)$ be a multiplicative Hom-Lie superalgebra. Then
the following statements hold:

$(1)$\quad ${\rm GDer}(L),{\rm QGer}(L)$ and ${\rm C}(L)$ are
Hom-subalgebras of  $\mho$.

$(2)$\quad ${\rm ZDer}(L)$ is a Hom-ideal of ${\rm Der}(L)$.}

 {\it Proof.}\quad  $(1)$ Assume that
$D_{\theta}\in{\rm GDer}_{\alpha^{k}}(L),~D_{\mu}\in{\rm
GDer}_{\alpha^{s}}(L),~\forall x\in {\rm hg}(L)$ and $y\in L.$ We
have
$$\begin{array}{ll}[(\tilde \alpha(D_{\theta}))(x),\alpha^{k+1}(y)]&=[(D_{\theta}\alpha)(x),\alpha^{k+1}(y)]=\alpha[D_{\theta}(x),\alpha^{k}(y)]\\
&=\alpha(D''_{\theta}([x,y])-(-1)^{\theta
x}[\alpha^{k}(x),D'_{\theta}(y)])\\&=\tilde
\alpha(D_{\theta}^{''})([x,y])-(-1)^{\theta
x}[\alpha^{k+1}(x),\tilde \alpha(D_{\theta}^{'})(y)]).\end{array}$$

Since both $\tilde \alpha(D_{\theta}^{''})$ and $\tilde
\alpha(D_{\theta}^{'})$ are in ${\rm End}_{\theta}(L),~\tilde
\alpha(D_{\theta})\in {\rm GDer}_{\alpha^{k+1}}(L)$ of degree
$\theta.$

We also
have$$\begin{array}{ll}[D_{\theta}D_{\mu}(x),\alpha^{k+s}(y)]&=D^{''}_{\theta}D^{''}_{\mu}([x,y])+(-1)^{\theta(\mu+x)}(-1)^{\mu
x}[\alpha^{s+k}(x),D^{'}_{\mu}D^{'}_{\theta}(y)]\\&-(-1)^{\theta(\mu+x)}D^{''}_{\mu}([\alpha^{k}(x),D^{'}_{\theta}(y)])-(-1)^{\mu
x}D^{''}_{\theta}([\alpha^{s}(x),D^{'}_{\mu}(y)])\end{array}$$
and$$\begin{array}{ll}[D_{\mu}D_{\theta}(x),\alpha^{k+s}(y)]&=D^{''}_{\mu}D^{''}_{\theta}([x,y])+(-1)^{\mu(\theta+x)}(-1)^{(\theta
x)}[\alpha^{s+k}(x),D^{'}_{\theta}D^{'}_{\mu}(y)]\\&-(-1)^{\mu(\theta+x)}D^{''}_{\theta}([\alpha^{s}(x),D^{'}_{\mu}(y)])-(-1)^{\theta
x}D^{''}_{\mu}([\alpha^{k}(x),D^{'}_{\theta}(y)]).\end{array}$$

\noindent Thus for all $x\in {\rm hg}(L)$ and $y\in L$, we
have$$[[D_{\theta},D_{\mu}](x),\alpha^{k+s}(y)]=[D''_{\theta},D''_{\mu}]([x,y])-(-1)^{(\theta+\mu)x}[\alpha^{k+s}(x),[D'_{\theta},D'_{\mu}](y)].$$

Since both $[D'_{\theta},D'_{\mu}]$ and $[D''_{\theta},D''_{\mu}]$
are in ${\rm End}_{\theta+\mu}(L)$, $[D_{\theta},D_{\mu}]\in {\rm
GDer}_{\alpha^{k+s}}(L)$ of degree $\theta+\mu$, $\forall \theta,\mu
\in G$, ${\rm GDer}(L)$ is a Hom-subalgebra of $ \mho$.

Similarly, ${\rm QGer}(L)$ is a Hom-subalgebra of $\mho$.

Assume that $D_{\theta}\in{\rm C}_{\alpha^{k}}(L),D_{\mu}\in{\rm
C}_{\alpha^{s}}(L),\forall x\in {\rm hg}(L),~y\in L.$ Then
$$\begin{array}{ll}[\tilde \alpha(D_{\theta})(x),\alpha^{k+1}(y)]&=\alpha([D_{\theta}(x),\alpha^{k}(y)])=(-1)^{\theta x}\alpha([\alpha^{k}(x),D_{\theta}(y)])\\
&=(-1)^{\theta x}[\alpha^{k+1}(x),\tilde \alpha
(D_{\theta})(y)],\end{array}$$ and so$$\tilde \alpha (D_{\theta})\in
{\rm C}_{\alpha^{k+1}}(L).$$
Note that $$\begin{array}{ll}[[D_{\theta},D_{\mu}](x),\alpha^{k+s}(y)]&=[D_{\theta}D_{\mu}(x),\alpha^{k+s}(y)]-(-1)^{\theta \mu}[D_{\mu}D_{\theta}(x),\alpha^{k+s}(y)]\\
&=D_{\theta}([D_{\mu}(x),\alpha^{s}(y)])-(-1)^{\theta \mu}D_{\mu}([D_{\theta}(x),\alpha^{k}(y)])\\
&=D_{\theta}D_{\mu}([x,y])-(-1)^{\theta
\mu}D_{\mu}D_{\theta}([x,y])\\&=[D_{\theta},D_{\mu}]([x,y]).\end{array}$$
Similarly,
$$(-1)^{(\theta+\mu)x}[\alpha^{k+s}(x),[D_{\theta},D_{\mu}](y)]=[D_{\theta},D_{\mu}]([x,y]).$$
Then $[D_{\theta},D_{\mu}]\in {\rm C}_{\alpha^{k+s}}(L)$ of degree
$\theta+\mu$, $\forall~\theta,\mu \in G.$ Thus ${\rm C}(L)$ is a
Hom-subalgebra of $\mho$.

$(2)$ Assume that $D_{\theta}\in{\rm
ZDer}_{\alpha^{k}}(L),D_{\mu}\in{\rm Der}_{\alpha^{s}}(L),\forall
x\in {\rm hg}(L),~y\in L.$ Then
$$\begin{array}{ll}[\tilde \alpha(D_{\theta})(x),\alpha^{k+1}(y)]&=\alpha([D_{\theta}(x),\alpha^{k}(y)])=\alpha D_{\theta}([x,y])=\tilde \alpha(D_{\theta})([x,y])=0.\end{array}$$ So$$\tilde \alpha (D_{\theta})\in
{\rm ZDer}_{\alpha^{k+1}}(L).$$
Note that $$\begin{array}{ll}[[D_{\theta},D_{\mu}]([x,y])]&=D_{\theta}D_{\mu}([x,y])-(-1)^{\theta \mu}D_{\mu}D_{\theta}([x,y])\\
&=D_{\theta}([D_{\mu}(x),\alpha^{s}(y)]+(-1)^{\mu
x}[\alpha^{s}(x),D_{\mu}(x)])\\&=0\end{array}$$
and$$\begin{array}{ll}[[D_{\theta},D_{\mu}](x),\alpha^{s+k}(y)]&=[(D_{\theta}D_{\mu}-(-1)^{\theta
\mu}D_{\mu}D_{\theta})(x),\alpha^{s+k}(y)]\\&=-(-1)^{\theta
\mu}(D_{\mu}([D_{\theta}(x),\alpha^{k}y]-(-1)^{\mu(\theta+x)}[\alpha^{s}(D_{\theta}(x)),D_{\mu}(\alpha^{k}y)])\\&=-(-1)^{\mu(\theta+x)}[D_{\theta}(\alpha^{s}(x)),\alpha^{k}(D_{\mu}(y))]=0.\end{array}$$
Then $[D_{\theta},D_{\mu}]\in {\rm ZDer}_{\alpha^{k+s}}(L)$ of
degree $\theta+\mu$, $\forall \theta,\mu \in G.$ Thus ${\rm
ZDer}(L)$ is a Hom-ideal of ${\rm Der}(L)$.\hfill$\Box$\vspace{0.3cm}

\noindent{\bf Lemma 2.2}\quad {\it Let  $(L,[\cdot,\cdot],\alpha)$
be a multiplicative Hom-Lie superalgebra. Then

    $(1)$ \quad $[{\rm Der}(L),{\rm C}(L)]\subseteq {\rm C}(L).$

    $(2)$ \quad $[{\rm QDer}(L),{\rm QC}(L)]\subseteq {\rm QC}(L).$

    $(3)$ \quad $[{\rm QC}(L),{\rm QC}(L)]\subseteq {\rm QDer}(L).$

    $(4)$ \quad ${\rm C}(L)\subseteq {\rm QDer}(L).$}

   {\it Proof.}\quad $(1)$\quad Assume that
$D_{\theta}\in{\rm GDer}_{\alpha^{k}}(L),D_{\mu}\in{\rm
C}_{\alpha^{s}}(L),\forall x\in {\rm hg}(L)$ and $y\in L.$ We have
$$\begin{array}{ll}[D_{\theta}D_{\mu}(x),\alpha^{k+s}(y)]&=D_{\theta}([D_{\mu}(x),\alpha^{s}(y)])-(-1)^{\theta(\mu+x)}[\alpha^{k}(D_{\mu}(x)),D_{\theta}(\alpha^{s}(y))]\\
&=D_{\theta}([D_{\mu}(x),\alpha^{s}(y)])-(-1)^{\theta(\mu+x)}[D_{\mu}(\alpha^{k}(x)),\alpha^{s}(D_{\theta}(y))]\\&=D_{\theta}D_{\mu}([x,y])-(-1)^{\theta(\mu+x)}(-1)^{\mu
x}[\alpha^{k+s}(x),D_{\mu}D_{\theta}(y)],\end{array}$$ and
$$\begin{array}{ll}[D_{\mu}D_{\theta}(x),\alpha^{k+s}(y)]&=D_{\mu}([D_{\theta}(x),\alpha^{k}(y)])\\
&=D_{\mu}D_{\theta}([x,y])-(-1)^{\theta x}D_{\mu}([\alpha^{k}(x),D_{\theta}(y)])\\
&=D_{\mu}D_{\theta}([x,y])-(-1)^{(\theta+\mu)x}[\alpha^{k+s}(x),D_{\mu}D_{\theta}(y)].\end{array}$$
Hence,
\begin{align*}[[D_{\theta},D_{\mu}](x),\alpha^{k+s}(y)]&=[D_{\theta}D_{\mu}(x),\alpha^{k+s}(y)]-(-1)^{\theta \mu}[D_{\mu}D_{\theta}(x),\alpha^{k+s}(y)]\\
&=D_{\theta}D_{\mu}([x,y])-(-1)^{\theta \mu}D_{\mu}D_{\theta}([x,y])\\
&=[D_{\theta},D_{\mu}]([x,y]).\end{align*} On the other hand,
\begin{align*}[D_{\theta}D_{\mu}(x),\alpha^{k+s}(y)]&=D_{\theta}([D_{\mu}(x),\alpha^{s}(y)])-(-1)^{\theta(\mu+x)}[\alpha^{k}(D_{\mu}(x)),D_{\theta}(\alpha^{s}(y))]\\
&=(-1)^{\mu x}(D_{\theta}([\alpha^{s}(x),D_{\mu}(y)])-(-1)^{\theta(\mu+x)}[\alpha^{k+s}(x),D_{\mu}D_{\theta}(y)])\\
&=(-1)^{\mu x}[D_{\theta}(\alpha^{s}(x)),\alpha^{k}(D_{\mu}(y))]+(-1)^{(\theta+\mu)x}[\alpha^{k+s}(x),D_{\theta}D_{\mu}(y)]\\
&\quad -(-1)^{\theta\mu}(-1)^{(\theta+\mu)x}[\alpha^{k+s}(x),D_{\mu}D_{\theta}(y)],\end{align*}
$$[D_{\mu}D_{\theta}(x),\alpha^{k+s}(y)]=(-1)^{\mu(\theta+x)}[D_{\theta}(\alpha^{s}(x)),\alpha^{k}(D_{\mu}(y))].$$
Then
\begin{align*}[[D_{\theta},D_{\mu}](x),\alpha^{k+s}(y)]&=[D_{\theta}D_{\mu}(x),\alpha^{k+s}(y)]-(-1)^{\theta \mu}[D_{\mu}D_{\theta}(x),\alpha^{k+s}(y)]\\
&=(-1)^{(\theta+\mu)x}([\alpha^{k+s}(x),D_{\theta}D_{\mu}(y)]-(-1)^{\theta \mu}[\alpha^{k+s}(x),D_{\mu}D_{\theta}(y)])\\
&=(-1)^{(\theta+\mu)x}[\alpha^{k+s}(x),[D_{\theta},D_{\mu}](y)].\end{align*}
Thus $[D_{\theta},D_{\mu}]\in{\rm C}_{\alpha^{s+k}}(L)$ of degree
$\theta+\mu,$ and we get $[{\rm Der}(L),{\rm C}(L)]\subseteq {\rm
C}(L).$

$(2)$\quad Similar to the proof of $(1)$.

$(3)$\quad  Assume that $D_{\theta}\in{\rm
QC}_{\alpha^{k}}(L),D_{\mu}\in{\rm C}_{\alpha^{s}}(L),\forall x\in
{\rm hg}(L)$ and $y\in L.$
 By Proposition~ 5.3~in~\cite{2}, we
have$$[[D_{\theta},D_{\mu}](x),\alpha^{k+s}(y)]+(-1)^{(\theta+\mu)x}[\alpha^{k+s}(x),[D_{\theta},D_{\mu}](y)]=0.$$
Let $D^{'}=0,~$hence $[D_{\theta},D_{\mu}]\in {\rm
QDer}_{\alpha^{k+s}}(L)$ of degree $\theta+\mu$ as desired.

$(4)$\quad See Proposition~5.2~in~\cite{2}.\hfill$\Box$
 \vspace{0.3cm}

\noindent{\bf Theorem 2.3}\quad {\it Let  $(L,[\cdot,\cdot],\alpha)$
be a multiplicative Hom-Lie superalgebra. Then $${\rm GDer}(L)={\rm
QDer}(L)+{\rm QC}(L).$$}

{\it Proof.}\quad Let $D_{\theta}\in {\rm GDer}_{\alpha^{k}}(L).$
 Then for all $x,y\in {\rm hg}(L)$, there exist
$D'_{\theta},D''_{\theta}\in {\rm End}_{\theta}(L)$ such that
$$[D_{\theta}(x),\alpha^{k}(y)]+(-1)^{\theta x}[\alpha^{k}(x),D'_{\theta}(y)]=D''_{\theta}([x,y]).$$
Since $(-1)^{(\theta+x)y}[\alpha^{k}(y),D_{\theta}(x)]+(-1)^{xy}[D'_{\theta}(y),\alpha^{k}(x)]=(-1)^{xy}D''_{\theta}([y,x]),$
$$[D'_{\theta}(y),\alpha^{k}(x)]+(-1)^{\theta y}[\alpha^{k}(y),D_{\theta}(x)]=D''_{\theta}([y,x]).$$
Hence $D'_{\theta}\in {\rm GDer}_{\alpha^{k}}(L).$ For all $x,y\in
{\rm hg}(L)$, we have
$$[\frac{D_{\theta}+D'_{\theta}}{2}(x),\alpha^{k}(y)]+(-1)^{\theta x}[\alpha^{k}(x),\frac{D_{\theta}+D'_{\theta}}{2}(y)]=D''_{\theta}([x,y]),$$
and
$$[\frac{D_{\theta}-D'_{\theta}}{2}(x),\alpha^{k}(y)]-(-1)^{\theta x}[\alpha^{k}(x),\frac{D_{\theta}-D'_{\theta}}{2}(y)]=0,$$
which imply that $\frac{D_{\theta}+D'_{\theta}}{2}\in {\rm
QDer}_{\alpha^{k}}(L)$ and $\frac{D_{\theta}-D'_{\theta}}{2}\in {\rm
QC}_{\alpha^{k}}(L),~of~degree~\theta$. Hence$$D_{\theta}=\frac{D_{\theta}+D'_{\theta}}{2}+\frac{D_{\theta}-D'_{\theta}}{2}\in
{\rm QDer}(L)+{\rm QC}(L),$$
 and$${\rm GDer}(L)\subseteq{\rm QDer}(L)+{\rm QC}(L).$$
 It is easy to vertify that ${\rm QDer}(L)+{\rm QC}(L)\subseteq {\rm GDer}(L).$
 Therefore ${\rm QDer}(L)+{\rm QC}(L)={\rm GDer}(L).$\hfill$\Box$
 \vspace{0.3cm}

\noindent{\bf Theorem 2.4}\quad {\it Let  $(L,[\cdot,\cdot],\alpha)$ be
a multiplicative Hom-Lie superalgebra, $\alpha$ a surjection and
${\rm Z}(L)$ the center of $L$. Then $[{\rm C}(L),{\rm
QC}(L)]\subseteq {\rm End}(L,{\rm Z}(L)).$
 Moreover, if ${\rm
Z}(L)=\{0\},$ then $[{\rm C}(L),{\rm QC}(L)]=\{0\}.$}

{\it Proof.}\quad Assume that $D_{\theta}\in {\rm
C}_{\alpha^{k}}(L), D_{\mu}\in {\rm QC}_{\alpha^{s}}(L)$ and $x\in
{\rm hg}(L)$.~Since $\alpha$ is surjection, $\forall y^{'}\in
L,~\exists y\in L,$ such that $y^{'}=\alpha^{k+s}(y).$  Then
$$\begin{array}{ll}[[D_{\theta},D_{\mu}](x),y^{'}]&=[[D_{\theta},D_{\mu}](x),\alpha^{k+s}(y)]\\&=[D_{\theta}D_{\mu}(x),\alpha^{k+s}(y)]-(-1)^{\theta\mu}[D_{\mu}D_{\theta}(x),\alpha^{k+s}(y)]\\
&=D_{\theta}([D_{\mu}(x),\alpha^{s}(y)])-(-1)^{\mu x}[\alpha^{s}D_{\theta}(x),D_{\mu}\alpha^{k}(y)]\\
&=D_{\theta}([D_{\mu}(x),\alpha^{s}(y)])-(-1)^{\mu x}D_{\theta}([\alpha^{s}(x),D_{\mu}(y)])\\
&=D_{\theta}([D_{\mu}(x),\alpha^{s}(y)]-(-1)^{\mu
x}[\alpha^{s}(x),D_{\mu}(y)])=0.\end{array}$$ Hence
$[D_{\theta},D_{\mu}](x)\in {\rm Z}(L)$, and
$[D_{\theta},D_{\mu}]\in {\rm End}(L,{\rm Z}(L))$ as desired.
Furthermore, if ${\rm Z}(L)=\{0\},$ it is clear that $[{\rm
C}(L),{\rm QC}(L)]=\{0\}.$ \hfill$\Box$
 \vspace{0.3cm}

By Theorem~2.3 in \cite{3}, if $(L,[\cdot,\cdot])$ is a Lie
superalgebra with $~{\rm
Z}(L)=\{0\},$ where ${\rm Z}(L)$ is the center of $L$, then ${\rm C}(L)={\rm QDer}(L)\cap{\rm QC}(L).$ But
it is not true in case that   $(L,[\cdot,\cdot],\alpha)$ is a
multiplicative Hom-Lie superalgebra.

\noindent{\bf  Example 2.5}\quad  Let $\{x_{1},x_{2},x_{3}\}$ be a
basis of a 3-dimensional linear space $L$ over $\bf K$. The following
bracket and linear map $\alpha$ on $L$ define a Hom-Lie algebra over
$\bf K$:
$$[x_{1},x_{2}]=~x_{1},~~~~~\alpha(x_{1})=~x_{1},$$
$$[x_{1},x_{3}]=~x_{2},~~~~~\alpha(x_{2})=2x_{2},$$
$$~[x_{2},x_{3}]=2x_{3},~~~~~~\alpha(x_{3})=2x_{3},$$
with $[x_{2},x_{1}],~[x_{3},x_{1}],~[x_{3},x_{2}]$ defined via
skewsymmetry.

Define $D:L\rightarrow L$ satisfying
$$D(x_{1})=x_{1},~D(x_{2})=2^{k}x_{2},~D(x_{3})=2^{k}x_{3}.~~(k\in Z_+)$$

It is obvious that ${\rm Z}(L)=0.$ $\forall y\in L,$ suppose
$y=ax_{1}+bx_{2}+cx_{3}.$ Define $ D^{'}\in {\rm End}(L)$ by
$$D^{'}(x_{1})=2^{k+1}x_{1},~D^{'}(x_{2})=2^{k+1
}x_{2},~D^{'}(x_{3})=2^{k+1}2^{k}x_{3}.$$

\noindent It is obvious that for
$i=1,2,3,$$$D^{'}([x_{i},ax_{1}+bx_{2}+cx_{3}])=[D(x_{i}),\alpha^{k}(ax_{1}+bx_{2}+cx_{3})]+[\alpha^{k}(x_{i}),D(ax_{1}+bx_{2}+cx_{3})],$$
and
$$[D(x_{i}),\alpha^{k}(ax_{1}+bx_{2}+cx_{3})]=[\alpha^{k}(x_{i}),D(ax_{1}+bx_{2}+cx_{3})].$$
But for all $t\in \bf Z$, $$D\in {\rm QDer}(L)\cap{\rm QC}(L).$$
 While$$D([x_{1},ax_{1}+bx_{2}+cx_{3}])=D(bx_{1}+cx_{2})=bx_{1}+c2^{k}x_{2}.$$
 $$[\alpha^{t}(x_{1}),D(ax_{1}+bx_{2}+cx_{3})]=[x_{1},ax_{1}+b2^{k}x_{2}+c2^{k}x_{3}]]=b2^{k}x_{1}+c2^{k}x_{2}.$$
That means $D\not\in {\rm C}(L).$

 \vspace{0.3cm}
\noindent{\bf  Proposition 2.6}\quad {\it Let
$(L,[\cdot,\cdot],\alpha)$ be a multiplicative Hom-Lie superalgebra and
$\alpha$ a surjection. If $~{\rm Z}(L)=\{0\}$, then ${\rm QC}(L)$
is a Hom-Lie superalgebra if and only if $~[{\rm QC}(L),{\rm
QC}(L)]=0.$}

{\it Proof.}\quad $(\Rightarrow)$\quad Assume that $D_{\theta}\in
{\rm hg}({\rm QC}_{\alpha^{k}}(L)),~D_{\mu}\in {\rm hg}({\rm
QC}_{\alpha^{s}}(L)),~x\in {\rm hg}(L).$ Since $\alpha$ is a
surjection, $\forall y^{'}\in L,~\exists y\in L,$ such that
$y^{'}=\alpha^{k+s}(y).$ Since ${\rm QC}(L)$ is a Hom-Lie
superalgebra, $[D_{\theta},D_{\mu}]\in {\rm hg}({\rm
QC}_{\alpha^{k+s}}(L))$, of degree $\theta+\mu.$ Then
$$[[D_{\theta},D_{\mu}](x),y^{'}]=[[D_{\theta},D_{\mu}](x),\alpha^{k+s}(y)]=(-1)^{(\theta+\mu)x}[\alpha^{k+s}(x),[D_{\theta},D_{\mu}](y)].$$
From the proof of Lemma~2.2~(1), we have
$$[[D_{\theta},D_{\mu}](x),y^{'}]=[[D_{\theta},D_{\mu}](x),\alpha^{k+s}(y)]=-(-1)^{(\theta+\mu)x}[\alpha^{k+s}(x),[D_{\theta},D_{\mu}](y)].$$
Hence
$[[D_{\theta},D_{\mu}](x),y^{'}]=[[D_{\theta},D_{\mu}](x),\alpha^{k+s}(y)]=0,$
i.e. $[D_{\theta},D_{\mu}]=0.$

$(\Leftarrow)$\quad It is clear.  \hfill$\Box$
\vspace{0.3cm}

\noindent{\bf  Definition 2.7} \,\cite{2}\quad {\it Let
$(L,\mu,\alpha)$ be a Hom-superalgebra.

$(1)$ The Hom-associator of $L$ is the trilinear map
$as_{\alpha}:L\times L\times L\rightarrow L$ defined
as$$as_{\alpha}=\mu\circ(\mu\otimes\alpha-\alpha\otimes\mu).$$ In
terms of elements, the map $as_{\alpha}$ is given by
$$as_{\alpha}(x,y,z)=\mu(\mu(x,y),\alpha(z))-\mu(\alpha(x),\mu(y,z))$$
for all $x,y,z\in hg(L).$

$(2)$ Let $L$ be a Hom-algebra over a field $\bf K$ of
characteristic $\not=2$ with an even bilinear multiplication
$\circ.$ If $L$ is $Z_{2}$-graded and $\alpha:L\rightarrow L$ is an
even linear map, then $(L,\circ,\alpha)$ is a Hom-Jordan
superalgebra if the following identities$$x\circ y=(-1)^{xy}y\circ x,$$
$$(-1)^{z(x+w)}as_{\alpha}(x\circ y,\alpha(z),\alpha(w))+(-1)^{x(y+z)}as_{\alpha}(y\circ w,\alpha(z),\alpha(x))$$
$$+(-1)^{y(w+z)}as_{\alpha}(w\circ
x,\alpha(z),\alpha(y))=0,$$hold for all $x,y,z\in hg(L).$}
\vspace{0.3cm}

\noindent{\bf  Proposition 2.8} \cite{2}\quad $(1)$ {\it Let
$(L,[\cdot,\cdot],\alpha)$ be a multiplicative Hom-Lie superalgebra,
with the operation $D_{\lambda}\bullet
D_{\theta}=D_{\lambda}D_{\theta}+(-1)^{\lambda\theta}D_{\theta}D_{\lambda},$
for all $\alpha$-derivations $D_{\lambda},D_{\theta}\in {\rm
hg}({\rm End}(L)),$ the triple $({\rm End}(L),\bullet,\alpha)$ is a
Hom-Jordan superalgebra.

 $(2)$
$(L,[\cdot,\cdot],\alpha)$ be a multiplicative Hom-Lie superalgebra,
with the operation $D_{\lambda}\bullet
D_{\theta}=D_{\lambda}D_{\theta}+(-1)^{\lambda\theta}D_{\theta}D_{\lambda},~$
for all elements $D_{\lambda},D_{\theta}\in {\rm hg}({\rm QC}(L))$.
 Then the triple $({\rm QC}(L),\bullet,\alpha)$ is a Hom-Jordan
superalgebra.} \vspace{0.3cm}

\noindent{\bf Theorem 2.9}\quad  {\it Let $(L,[\cdot,\cdot],\alpha)$
be a multiplicative Hom-Lie superalgebra. Then ${\rm QC}(L)$ is a
Hom-Lie superalgebra with $[D_{\lambda},D_{\theta}]=
D_{\lambda}D_{\theta}-(-1)^{\lambda\theta}D_{\theta}D_{\lambda}$ if
and only if ${\rm QC}(L)$ is also a Hom-associative superalgebra
with respect to composition.}

{\it Proof.}\quad $(\Leftarrow)$ For all $D_{\lambda},D_{\theta}\in
{\rm hg}({\rm QC}(L)),$ we have $D_{\lambda}D_{\theta}\in {\rm
QC}(L)$ and $D_{\theta}D_{\lambda}\in {\rm QC}(L)$, so
$[D_{\lambda},D_{\theta}]=D_{\lambda}D_{\theta}-(-1)^{\lambda\theta}D_{\theta}D_{\lambda}\in
{\rm QC}(L).$ Hence, ${\rm QC}(L)$ is a Hom-Lie superalgebra.

$(\Rightarrow)$ Note that $D_{\lambda}D_{\theta}=D_{\lambda}\bullet
D_{\theta}+\frac{[D_{\lambda},D_{\theta}]}{2}$ and by Proposition
2.8, we have $D_{\lambda}\bullet D_{\theta}\in {\rm QC}(L),~
[D_{\lambda},D_{\theta}]\in QC(L)$. It follows that
$D_{\lambda}D_{\theta}\in {\rm QC}(L)$ as desired.\hfill$\Box$
    \vspace{0.3cm}

\noindent {\bf \S 3\quad  The quasiderivations of Hom-Lie
superalgebras }
        \vspace{0.3cm}

       In this section, we will prove that the quasiderivations of $L$ can
be embedded as derivations in a larger Hom-Lie superalgebra and
obtain a direct sum decomposition of Der($L$) when the annihilator
of $L$ is equal to zero.

 \noindent{\bf Proposition 3.1}\quad {\it Let $(L,[\cdot,\cdot],\alpha)$ be a Hom-Lie superalgebra over
       ${\bf K}$ and $t$ an indeterminate. We define $\breve{L}_{g}:=L_{g}[tF[t]/(t^{3})]=
\{\Sigma(x_{g}\otimes t+y_{g}\otimes t^{2})| x_{g},y_{g}\in
L_{g}\},~\breve\alpha(\breve{L}_{g}):=\{\Sigma(\alpha(x_{g})\otimes
t+\alpha(y_{g})\otimes t^{2}):x_{g},y_{g}\in L_{g}\}$, and let
$\breve{L}=\breve{L}_{\bar0}\oplus\breve{L}_{\bar1}$. Then
$\breve{L}$ is a Hom-Lie superalgebra with the operation
$[x_{\lambda}\otimes t^{i},x_{\theta}\otimes
t^{j}]=[x_{\lambda},x_{\theta}]\otimes t^{i+j},$ for all
$x_{\lambda},x_{\theta}\in {\rm hg}(L),i,j\in\{1,2\}$}.

{\it Proof.}\quad For all $x_{\lambda},x_{\theta},x_{\mu}\in {\rm
hg}(L)$ and $i,j,k \in\{1,2\},$ we have
$$\begin{array}{ll}[x_{\lambda}\otimes t^{i},x_{\theta}\otimes t^{j}]&=[x_{\lambda},x_{\theta}]\otimes t^{i+j}\\
&=-(-1)^{\lambda\theta}[x_{\theta},x_{\lambda}]\otimes t^{i+j}\\
&=-(-1)^{\lambda\theta}[x_{\theta}\otimes t^{j},x_{\lambda}\otimes
t^{i}],\end{array}$$ and \begin{eqnarray*}
&&\ \ [\breve\alpha(x_{\lambda}\otimes t^{i}),[x_{\theta}\otimes t^{j},x_{\mu}\otimes t^{k}]]=[\alpha(x_{\lambda}),[x_{\theta},x_{\mu}]]\otimes t^{i+j+k}\\
&&=([[x_{\lambda},x_{\theta}],\alpha(x_{\mu})]+(-1)^{\lambda\theta}[\alpha(x_{\theta}),[x_{\lambda},x_{\mu}]])\otimes t^{i+j+k}\\
&&=[[x_{\lambda},x_{\theta}],\alpha(x_{\mu})]\otimes t^{i+j+k}+(-1)^{\lambda\theta}[\alpha(x_{\theta}),[x_{\lambda},x_{\mu}]]\otimes t^{i+j+k}\\
&&=[[x_{\lambda}\otimes t^{i},x_{\theta}\otimes
t^{j}],\breve\alpha(x_{\mu}\otimes t^{k})]+
(-1)^{\lambda\theta}[\breve\alpha(x_{\theta}\otimes
t^{j}),[x_{\lambda}\otimes t^{i},x_{\mu}\otimes
t^{k}]].\end{eqnarray*} Hence $\breve{L}$ is a Hom-Lie superalgebra.
\hfill$\Box$ \vspace{0.3cm}

For notational convenience, we write $xt(xt^{2})$ in place of $x\otimes t(x\otimes t^{2}).$

If $U$ is a $\rm Z_{2}$-graded subspace of $L$ such that $L=U\oplus
[L,L],$ then $$\breve{L}=Lt+Lt^{2}=Lt+[L,L]t^{2}+Ut^{2},$$ or more
precisely,$$\breve{L}=\breve{L}_{\bar0}\oplus\breve{L}_{\bar1}
=(L_{\bar0}t+[L,L]_{\bar0}t^{2}+U_{\bar0}t^{2})\oplus(L_{\bar1}t+[L,L]_{\bar1}t^{2}+U_{\bar1}t^{2}).$$

Now we define a map $\varphi:{\rm QDer}(L)\rightarrow {\rm
End}(\breve{L})$ satisfying
$$\varphi(D)(at+bt^{2}+ut^{2})=D(a)t+D'(b)t^{2},$$
where $D, D'$ satisfy
$(1.2)$, $a\in {\rm hg}(L),b\in {\rm hg}([L,L]),u\in {\rm hg}(U)$
and $d(a)=d(b)=d(u).$ \vspace{0.3cm}

 \noindent{\bf Proposition 3.2}\quad {\it
$(1)$\quad $d(\varphi)=0.$

$(2)$\quad $\varphi$ is injective and $\varphi(D)$ does not depend on the choice of $D'$.

$(3)$\quad $\varphi({\rm QDer}(L))\subseteq {\rm Der}(\breve{L}).$}

{\it Proof.}\quad (1)\quad It is clear.

(2)\quad If $\varphi(D_{\lambda})=\varphi(D_{\theta}),$ then for all $a\in {\rm hg}(L),b\in {\rm hg}([L,L])$ and $u\in {\rm hg}(U),$ we have
$$\varphi(D_{\lambda})(at+bt^{2}+ut^{2})=\varphi(D_{\theta})(at+bt^{2}+ut^{2}),$$ that is $$D_{\lambda}(a)t+D'_{\lambda}(b)t^{2}=
D_{\theta}(a)t+D'_{\theta}(b)t^{2},$$ so
$D_{\lambda}(a)=D_{\theta}(a).$ Hence $D_{\lambda}=D_{\theta},$ and
$\varphi$ is injective.

Suppose that there exists $D''$ such that
$$\varphi(D)(at+bt^{2}+ut^{2})=D(a)t+D''(b)t^{2},$$ and
$$[D(x),\alpha^{k}(y)]+(-1)^{Dx}[\alpha^{k}(x),D(y)]=D''([x,y]),$$
then we have $$D'([x,y])=D''([x,y]),$$ thus $D'(b)=D''(b).$ Hence
$$\varphi(D)(at+bt^{2}+ut^{2})=D(a)t+D'(b)t^{2}=D(a)t+D''(b)t^{2},$$
which implies $\varphi(D)$ is determined by $D$.

(3)\quad We have $[x_{\lambda}t^{i},x_{\theta}t^{j}]=[x_{\lambda},x_{\theta}]t^{i+j}=0,$ for all $i+j\geq 3$.
Thus, to show $\varphi(D)\in {\rm Der}(\breve{L}),$ we need only to check the validness of the following equation
$$\varphi(D)([xt,yt])=[\varphi(D)(xt),\breve\alpha^{k}(yt)]+(-1)^{Dx}[\breve\alpha^{k}(xt),\varphi(D)(yt)].$$
For all $x,y\in{\rm hg}(L),$ we have
$$\begin{array}{ll}\varphi(D)([xt,yt])&=\varphi(D)([x,y]t^2)=D'([x,y])t^2\\
&=([D(x),\alpha^{k}(y)]+(-1)^{Dx}[\alpha^{k}(x),D(y)])t^2\\
&=[D(x)t,\alpha^{k}(y)t]+(-1)^{Dx}[\alpha^{k}(x)t,D(y)t]\\
&=[\varphi(D)(xt),\breve\alpha^{k}(yt)]+(-1)^{Dx}[\breve\alpha^{k}(xt),\varphi(D)(yt)].\end{array}$$
Therefore, for all $D\in {\rm QDer}(L)$, we have $\varphi(D)\in {\rm
Der}(\breve{L})$.\hfill$\Box$ \vspace{0.3cm}

\noindent{\bf Proposition 3.3}\quad {\it Let $(L,[\cdot,\cdot],\alpha)$
be a multiplicative Hom-Lie superalgebra and $\alpha$ a
surjection. ${\rm Z}(L)=\{0\}$ and $\breve{L},~\varphi$ are as
defined above. Then ${\rm Der}(\breve{L})=\varphi({\rm
QDer}(L))\oplus {\rm ZDer}(\breve{L}).$}

{\it Proof.}\quad Since ${\rm Z}(L)=\{0\}$, we have ${\rm
Z}(\breve{L})=Lt^2.$ For all $g\in {\rm Der}(\breve{L}),$ we have
$g({\rm Z}(\breve{L}))\subseteq {\rm Z}(\breve{L}),$ hence
$g(Ut^2)\subseteq g({\rm Z}(\breve{L}))\subseteq {\rm
Z}(\breve{L})=Lt^2.$ Now we define a map
$f:Lt+[L,L]t^2+Ut^2\rightarrow Lt^2$ by
$$\ f(x)=\left\{\begin{array}{ll}g(x)\cap Lt^2,& x\in Lt ;\\
 g(x),& x\in Ut^2 ;\\  0,& x\in [L,L]t^2.\end{array}\right.$$
It is clear that $f$ is linear. Note that
$$f([\breve{L},\breve{L}])=f([L,L]t^2)=0,~[f(\breve{L}),\breve\alpha^{k}{L}]\subseteq [Lt^2,\alpha^{k}(L)t+\alpha^{k}(L)t^2]=0,$$ hence $f\in {\rm ZDer}(\breve{L}).$
Since $$(g-f)(Lt)=g(Lt)-g(Lt)\cap Lt^{2}=g(Lt)-Lt^{2}\subseteq Lt,~
(g-f)(Ut^2)=0,$$ and
$$(g-f)([L,L]t^2)=g([\breve{L},\breve{L}])\subseteq
[\breve{L},\breve{L}]=[L,L]t^2,$$ there exist $D,~D'\in {\rm
End}(L)$ such that for all $a\in L,~b\in [L,L]$, $$(g-f)(at)=D(a)t,~
(g-f)(bt^2)=D'(b)t^2.$$ Since $(g-f)\in {\rm Der}(\breve{L})$ and by
the definition of ${\rm Der}(\breve{L})$, we have
$$[(g-f)(a_1t),\breve\alpha^{k}(a_2t)]+(-1)^{(g-f)a_1}[\breve\alpha^{k}(a_1t),(g-f)(a_2t)]=(g-f)([a_1t,a_2t]),$$
for all $a_1,a_2\in L.$ Hence
$$[D(a_1),\breve\alpha^{k}(a_2)]+(-1)^{Da_1}[\breve\alpha^{k}(a_1),D(a_2)]=D'([a_1,a_2]).$$ Thus
$D\in {\rm QDer}(L).$ Therefore, $$g-f=\varphi(D)\in \varphi({\rm
QDer}(L))\Rightarrow {\rm Der}(\breve{L})\subseteq \varphi({\rm
QDer}(L))+{\rm ZDer}(\breve{L}).$$ By Proposition $3.2~(3)$ we
have $${\rm Der}(\breve{L})=\varphi({\rm QDer}(L))+{\rm
ZDer}(\breve{L}).$$

For all $f\in \varphi({\rm QDer}(L))\cap{\rm ZDer}(\breve{L})$, there exists an element $D\in {\rm QDer}(L)$ such that $f=\varphi(D).$ Then
$$f(at+bt^2+ut^2)=\varphi(D)(at+bt^2+ut^2)=D(a)t+D'(b)t^2,$$ for all $a\in L,b\in [L,L].$

On the other hand, since $f\in {\rm ZDer}(\breve{L}),$ we have
$$f(at+bt^2+ut^2)\in {\rm Z}(\breve{L})=Lt^2.$$ That is to say,
$D(a)=0,$ for all $a\in L$ and so $D=0.$ Hence $f=0.$

Therefore ${\rm Der}(\breve{L})=\varphi({\rm QDer}(L))\oplus {\rm
ZDer}(\breve{L})$ as desired. \hfill$\Box$\vspace{4mm}
     \end{document}